\documentclass[final]{siamltex}
\usepackage{graphicx,amssymb,amsmath,mathrsfs}

\newtheorem{thm}{Theorem}[section]
\newtheorem{lem}{Lemma}[section]
\newtheorem{rem}{Remark}[section]

\author{
Gunther H. Peichl
\thanks{University of Graz, Institute for Mathematics, Heinrichstr. 36, 8010 Graz, Austria
({\tt gunther.peichl@kfunigraz.ac.at}).}
\and Rachid Touzani
\thanks{Laboratoire de Math\'ematiques, Universit\'e Blaise Pascal (Clermont-Ferrand) and
CNRS (UMR 6620), Campus Universitaire des C\'ezeaux, 63177 Aubi\`ere cedex, France.
({\tt Rachid.Touzani@univ-bpclermont.fr})}
}

\begin{document}
\title{An accurate finite element method for elliptic interface problems}



\maketitle

\begin{abstract}
A finite element method for elliptic problems with discontinuous
coefficients is presented. The discontinuity is assumed to take place along a closed smooth
curve. The proposed method allows to deal with meshes that are not adapted to
the discontinuity line. The (nonconforming) finite element space is enriched with local basis functions.
We prove an optimal convergence rate in the $H^1$--norm. Numerical tests confirm the theoretical results.
\end{abstract}

\section{Introduction}

Boundary value problems with discontinuous coefficients constitute a prototype of various problems
in heat transfer and continuum mechanics where heterogeneous media are involved. The numerical solution
of such problems requires much care since their solution does not generally enjoy enough smoothness properties
required to obtain optimal convergence rates. Although fitted or adapted meshes can handle such difficulties, these
solution strategies become expensive if the discontinuity front evolves with time or within an iterative process.
Such a (weak) singularity appears also in the numerical solution of other types of problems
like  free boundary problems when they are formulated for a fixed mesh or for fictitious domain methods.

We address, in this paper, a new finite element approximation of a model elliptic transmission problem that
allows nonfitted meshes. It is well known that the standard finite element approximation of such a problem does not
converge with a first order rate in the $H^1$-norm in the general case. We propose a method that converges optimally
provided the interface curve is a sufficiently smooth curve. Our method is based on a local enrichment of the finite
element space in the elements intersected by the interface. The local feature is ensured by the use of a hybrid approximation.
A Lagrange multiplier enables to recover the conformity of the approximation. The derived method appears then rather
as a local modification of the equations of interface elements than a modification of the linear system of equations.
This property ensures that the structure of the matrix of the linear system is not affected by the enrichment.

Let us mention other authors who addressed this topic in the finite element context. We point out the
so-called XFEM (eXtended Finite Element Methods) developed in Belytschko \emph{et al.} \cite{BMUP}
where the finite element space is  modified in interface elements by using the level set function
associated to the interface. Such methods, that are used also for crack propagation, have in our
point of view, the drawback of resulting in a variable matrix structure.
Moreover, although no theoretical analysis is available, numerical experiments
show that they are not optimal in terms of accuracy. Other authors like Hansbo {\it et al.} \cite{HH,HLPS},
have similar approaches to ours but here also the proposed method seems to
modify the matrix structure by enriching the finite element.
In Lamichhane--Wohlmuth \cite{LW} and Braess--Dahmen \cite{BD}, a similar Lagrange multiplier approach is used
for a mortar finite element
formulation of a domain decomposition method. Finally, in a work by Li \emph{et al}
\cite{LLW}, an immersed interface technique, inspired from finite difference schemes, is adapted to
the finite element context. Note also that the references where Lagrange multipliers are employed have
for these multipliers as supports the edges defining the interface. In our method, the interface supports the
added degrees of freedom but the Lagrange multipliers are defined on the edges intersected by the interface
and thus serve to compensate the nonconformity of the finite element space rather than enforcing
interface conditions, which are being naturally ensured by the variational formulation.

\bigskip\noindent
In the following, we  use the space $L^2(\Omega)$ equipped with the norm $\|\cdot\|_{0,\Omega}$
and the Sobolev spaces $H^m(\Omega)$ and $W^{m,p}(\Omega)$ endowed with the norms $\|\cdot\|_{m,\Omega}$ and
$\|\cdot\|_{m,p,\Omega}$ respectively. We shall also use the semi-norm $|\cdot|_{1,\Omega}$ of $H^1(\Omega)$.
Moreover, if $\Omega_1$ and $\Omega_2$ form a partition of $\Omega$, {\em i.e.}, $\overline\Omega=
\overline\Omega_1\cup\overline\Omega_2$, $\Omega_1\cap\Omega_2=\emptyset$ and if $v$ is
a function in $W^{m-1,p}(\Omega)$ with $v_{|\Omega_i}\in W^{m,p}(\Omega_i)$,
then we shall adopt the convention $v\in W^{m,p}(\Omega_1\cup\Omega_2)$ and denote
by $\|v\|_{m,p,\Omega_1\cup\Omega_2}$ the broken Sobolev norm
$$
\|v\|_{m,p,\Omega_1\cup\Omega_2} = \|v\|_{m-1,p,\Omega} +
\|v\|_{m,p,\Omega_1} + \|v\|_{m,p,\Omega_2}.
$$
Similarly, we denote by $\|\cdot\|_{m,\Omega_1\cup\Omega_2}$ and $|\cdot|_{m,\Omega_1\cup\Omega_2}$,
the broken Sobolev norm and semi-norm respectively for the $H^m$--space.
Finally, we shall denote by $C$, $C_1, C_2, \ldots$ various generic constants that
do not depend on mesh parameters and by $|A|$ the Lebesgue measure of a set $A$
and by $A^\circ$ the interior of a set $A$.

\bigskip\noindent
Let $\Omega$ denote a domain in $\mathbb R^2$ with smooth boundary $\Gamma$ and let $\gamma$ stand for
a closed $C^2$-curve in $\Omega$ which separates $\Omega$ into two disjoint subdomains $\Omega^+$,
$\Omega^-$ such that $\Omega= \Omega^+\cup\gamma\cup\Omega^-$ and $\partial \Omega^+ =\gamma$.
For given $f\in L^2(\Omega)$ and $a\in L^\infty(\Omega)$ we consider the transmission problem:
$$
\left\{
\begin{aligned}{}
-\nabla\cdot(a\nabla u) &= f&&\qquad\text{in }\Omega^+\cup\Omega^-,\\
u &= 0&&\qquad\text{on }\Gamma,\\
[u] = \big[a\frac{\partial u}{\partial n}\big] &= 0&&\qquad\text{on }\gamma,
\end{aligned}
\right.
$$
where $[v]$ denotes the jump of a quantity $v$ across the interface $\gamma$ and $n$ is the normal
unit vector to $\gamma$ pointing into $\Omega^-$. For definiteness we let $[v] = v^- -v^+$ with
$v^\pm = v_{|\Omega^\pm}$. In addition to boundedness of the diffusion coefficient we assume
\begin{align}\label{eq:ass-a}
\begin{split}
a^\pm &\in W^{1,\infty}(\Omega^\pm),\\
a(x) &\ge \alpha>0,\qquad\text{for } x\in\Omega,
\end{split}
\end{align}
i.e. $a$ is uniformly continuous on $\Omega\setminus \gamma$, but discontinuous across $\gamma$.

The standard variational formulation of this problem consists in seeking $u\in H^1_0(\Omega)$
such that
\begin{equation}
\int_\Omega a\,\nabla u\cdot\nabla v\,dx = \int_\Omega fv\,dx\qquad\forall\ v\in H^1_0(\Omega).
\label{Pb}
\end{equation}
In view of the ellipticity condition \eqref{eq:ass-a}, Problem \eqref{Pb} has a unique solution
$u$ in  $H^1_0(\Omega)$ but clearly $u\notin H^2(\Omega)$. We shall assume throughout this
paper the regularity properties:
\begin{align}
&u_{|\Omega^-}\in H^2(\Omega^-),\quad u_{|\Omega^+}\in H^2(\Omega^+),\nonumber\\
&\|u\|_{2,\Omega^-\cup\Omega^+} \le C\,\|f\|_{0,\Omega}.\label{H2-regularity}
\end{align}
Note that these assumptions are satisfied in the case where $a_{|\Omega^-}$ and
$a_{|\Omega^+}$ are constants (see \cite{Lemrabet,Petzoldt} for instance).

In the following, we describe a {\em fitted finite element method}.
defined by adding extra unknowns on the interface $\gamma$. It turns out that this method
leads to  an optimal convergence rate. Although it is well suited for the model problem
it  seems to be inefficient in more elaborate problems which, for example, involve moving interfaces.
To circumvent this difficulty, we define a new method where the added degrees of freedom have local
supports and then yield a nonconforming finite element method. We show that the use of a Lagrange
multiplier removes this nonconformity and ensures an optimal convergence rate.

\section{A fitted finite element method}

\begin{figure}[!ht]
\centering
\includegraphics*[height=4cm]{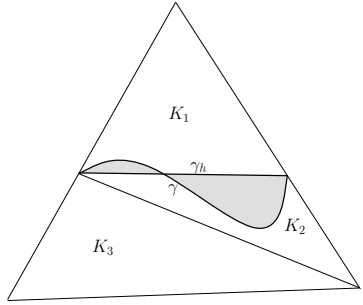}
\caption{\label{Fig-1} Subdivision of interface triangles.}
\end{figure}

Assume that the domain $\Omega$ is a convex polygon and consider a regular triangulation ${\mathscr T}_h$ of
$\overline\Omega$ with closed  triangles whose edges have lengths $\le h$. We assume that $h$ is  small enough so that for each triangle $T\in{\mathscr T}_h $ only the following cases have to be considered:
\begin{enumerate}
\item[1)] $T\cap \gamma = \emptyset$.
\item[2)] $T\cap \gamma$ is an edge or a vertex of $T$.
\item [3)]$\gamma$ intersects two different edges of $T$ in two distinct points different from the vertices.
\item[4)] $\gamma$ intersects one edge and its opposite vertex.
\end{enumerate}
 Let $V_h$ denote the lowest degree finite element
space
$$
V_h = \{\,v\in C^0(\overline\Omega);\ v_{|T}\in P_1(T)\ \forall\ T\in {\mathscr T}_h,\ v=0\text{ on }\Gamma\},
$$
where $P_1(T)$ is the space of affine functions on $T$. A finite element approximation of \eqref{Pb}
consists in computing $u_h\in V_h$ such that
\begin{equation}
\int_\Omega a\,\nabla u_h\cdot\nabla v\,dx = \int_\Omega fv\,dx\qquad\forall\ v\in V_h.
\label{SFEM}
\end{equation}
It is well known that, since $u\notin H^2(\Omega)$, the classical error estimates
(see \cite{Ciarlet}) do not hold any more even though we still have the convergence result,
$$
\lim_{h\to 0}\|u-u_h\|_{1,\Omega} = 0.
$$
A fitted treatment of the interface $\gamma$ can however improve this result. Let  for this purpose $\mathscr T_h^\gamma$
denote   the set of triangles that intersect the interface $\gamma$ corresponding to cases 3) and 4) above,
$$
\mathscr T_h^\gamma := \{T\in\mathscr T_h;\ \gamma\cap T^\circ\neq\emptyset\},
$$
and consider a continuous piecewise linear interpolation of $\gamma$, denoted by $\gamma_h$, as shown in Figure
\ref{Fig-1}. Clearly, $\gamma_h$ is the line that intersects $\gamma$ at two edges of any triangle that contains $\gamma$.
Unless the intersection of $\gamma$ with the boundary of a triangle $T$ does not coincide
with an edge, $T$ is split into two sets $T^+$ and $T^-$ separated
by the curve $\gamma$. In case 3), the straight line $\gamma_h\cap T$ splits $T$ into a triangle $K_1$
and a quadrilateral that we split into two subtriangles $K_2$ and $K_3$, where we choose $K_2$ such that
$K_1\cap K_2=\gamma_h$. In case 4), $\gamma_h\cap T$ splits $T$ into two triangles $K_1$ and $K_2$.
In this case we set $K_3=\emptyset$. This construction
 defines the new fitted finite element mesh of the domain $\Omega$ (see Figure \ref{Fig-1}).
The splitting $T=K_1\cup K_2\cup K_3$ is not unique but
the convergence analysis does not depend on it. Let us denote by $\mathscr T_T^\gamma$ the set of the
three subtriangles of $T$.  Below
$\mathscr E_h$ will stand for the set of all edges of elements and $\mathscr E_h^\gamma$ is the set of all edges
that are intersected by $\gamma$ (or $\gamma_h$), {\em i.e.}
$$
\mathscr E_h^\gamma := \{e\in\mathscr E_h;\ \gamma\cap e^\circ\neq\emptyset\}.
$$
For each $T\in\mathscr T_h$, $\mathscr E_T$ is the set of the three edges of $T$.
The fitted mesh is denoted by $\mathscr T_h^F$, {\em i.e.}
$$
\mathscr T_h^F := \mathscr T_h \cup \bigcup_{T\in\mathscr T_h^\gamma}\Big(\cup_{K\in\mathscr T_K^\gamma}K\Big).
$$
and by $S_h^\gamma := \bigcup\{T;\ T\in\mathscr T_h^\gamma\}$. Let us finally note that the curve $\gamma_h$
defines a new splitting of $\Omega$ into two subdomains $\Omega^-_h$ and $\Omega^+_h$
where $\Omega^\pm_h$ is defined analogously to $\Omega^\pm$ with $\gamma$ replaced by $\gamma_h$.

Next we  construct an approximation of the function $a$ on the elements of $\mathscr T_h^F$:
For this purpose, let $\tilde a^{\pm}$ be  extensions of $a^{\pm}$ to $\Omega$ such that
$\tilde a^{\pm}\in W^{1,\infty}(\Omega)$. Such extensions exist due to the regularity of $\gamma$
(see \cite{Adams}). Define $\tilde a_h\in W^{1,\infty}(\Omega)$ by
\[
\tilde a_h = \begin{cases}
\tilde a^+ & \text{in } \Omega_h^+,\\
\tilde a^- & \text{in } \Omega_h^-,
\end{cases}
\]
and denote by $a_h$ the piecewise linear interpolant of $\tilde a$ on $\mathscr T_h^F$. Hence $a_h$ is continuous on $\Omega_h^+\cup \Omega_h^-$ and coincides with $a$ on the nodes of $\mathscr T_h^F$.
In addition, the function $a_h$ is discontinuous
across the line $\gamma_h$ and satisfies the properties,
\begin{align}
&a_{h|\Omega^+_h}\in W^{1,\infty}(\Omega^+_h),\ a_{h|\Omega^-_h}\in W^{1,\infty}(\Omega^-_h),\label{Prop-ah-1}\\
&\|a_h\|_{0,\infty,\Omega} \le C\,\|a\|_{0,\infty,\Omega},\label{Prop-ah-2}\\
&a_h \ge\alpha > 0\qquad\text{a.e. in }\Omega.\label{Prop-ah-3}
\end{align}
We now define the finite element space
\begin{align*}
&W_h = V_h + X_h,\\
&X_h:=\{v\in C^0(\overline\Omega);\ v_{|\Omega\setminus S_h^\gamma}=0,\ v_{|K}\in
P_1(K)\ \forall\ K\in\mathscr T_T^\gamma,\ \forall\ T\in\mathscr T_h^\gamma\}.
\end{align*}
Note that we have $W_h\subset H^1_0(\Omega)$.
A fitted finite element approximation is defined as the follows:
\begin{equation}
\left\{
\begin{aligned}{}
&\text{Find }u_h^F\in W_h\text{ such that}\\
&\int_\Omega a_h\nabla u_h^F\cdot\nabla v\,dx = \int_\Omega fv\,dx\qquad\forall\ v\in W_h.
\end{aligned}
\right.
\label{AFEM}
\end{equation}

In order to study the convergence of Problem \eqref{AFEM}, we consider the auxiliary problem:
\begin{equation}
\left\{
\begin{aligned}{}
&\text{Find }\widehat u_h\in H^1_0(\Omega)\text{ such that}\\
&\int_\Omega a_h\nabla\widehat u_h\cdot\nabla v\,dx = \int_\Omega fv\,dx\qquad\forall\ v\in H^1_0(\Omega).
\end{aligned}
\right.
\label{Ah}
\end{equation}

We note that both problems \eqref{AFEM} as well as \eqref{Ah} have a unique solution.
The regularity properties \eqref{H2-regularity} imply $u^+\in C^0(\bar\Omega^+)$, $u^-\in C^0(\bar\Omega^-)$
and that $u^+$ and $u^-$ have a common trace on $\gamma$. Therefore $u$ is continuous on $\Omega$
and the piecewise $P_1$
interpolant $I_h u\in W_h$ is well defined.  In the following let $\tilde u^\pm \in H^2(\Omega)$ stand for the
extensions of $u^\pm$ from $\Omega^\pm$ to $\Omega$.

In the sequel, we assume that the fitted family of meshes $(\mathscr T_h\cup\mathscr T_h^\gamma)_h$
satisfies the condition
\begin{equation}
\frac h\varrho\le C\,h^{-\theta}\label{MeshRegularity}
\end{equation}
for some $\theta\in [0,1)$ and for which $C$ is independent of $h$, where $\varrho$ denotes the radius
of the largest ball contained in any triangle in any triangle $T\in\mathscr T_h^F$.

\begin{lem}\label{eq:interpolation-error}
Let $u\in H^2(\Omega^+\cup\Omega^-)$.
\begin{enumerate}
\item We have the local interpolation error
\begin{equation}
\label{InterpolationErrorLocal}
|u-I_h u|_{1,T} \le
\begin{cases}
Ch \,|u|_{2,T} & \text{for }T\in \mathscr T_h\setminus \mathscr T_h^\gamma\\
C \frac{h^2}{\varrho_K}(|\tilde u^+|_{2,K} + |\tilde u^-|_{2,K}) &\text{for }K\in \mathscr T_T^\gamma, \,T\in \mathscr T_h^\gamma,
\end{cases}
\end{equation}
where $\varrho_K$ is the radius of the inscribed circle of $K$.
\item The global interpolation error is given by
\begin{equation}
|u-I_hu|_{1,\Omega} \le C\,h^{1-\theta}\,|u|_{2,\Omega^+\cup\Omega^-}.\label{InterpolationErrorGlobal-1}
\end{equation}
Moreover, if $u\in W^{2,\infty}(\Omega^+\cup\Omega^-)$ then
\begin{equation}
|u-I_hu|_{1,\Omega} \le C\,h\,|u|_{2,\infty,\Omega^+\cup\Omega^-}.\label{InterpolationErrorGlobal-2}
\end{equation}
\end{enumerate}
\end{lem}

\begin{proof}
Since the local interpolation error estimate for $T\in\mathscr T_h\setminus\mathscr T_h^\gamma$ is
classic in finite element theory (see \cite{BS} or \cite{Ciarlet} for instance), we only need to prove
the second estimate on triangles where $u$ is only piecewise smooth. Consider an element $T\in \mathscr T_h^\gamma$ and any subtriangle $K\in \mathscr T_T^\gamma$. Without loss of generality we assume $K\subset \Omega_h^+$, then
\begin{align*}
K &= (K\cap \Omega^+) \,\cup\, (K\cap \Omega^-).
\end{align*}
Since $K\cap \Omega^-\subset T\cap\Omega^-\cap \Omega_h^+$ and $\gamma_h$ interpolates the interface $\gamma$ we obtain for the measure of $K\cap \Omega^-$
\begin{equation}\label{eq:measure-k-omega}
|K\cap \Omega^-| \le |T\cap\Omega^-\cap \Omega_h^+|\le  C h^3,
\end{equation}
with a constant $C>0$ which  depends on $\gamma$ only. In view of $I_h u = I_h \tilde u^+$, the standard
 interpolation theory (see \cite{Ciarlet} or \cite{BS}) implies
\begin{align}\label{eq:interp1}
\begin{split}
|u-I_h u|_{1,K} &\le |u-\tilde u^+|_{1,K} + |\tilde u^+-I_h \tilde u^+|_{1,K}\\
&\le
|u-\tilde u^+|_{1,K} + C\,\frac{h^2}{\varrho_K}\,|\tilde u^+|_{2,K}.
\end{split}
\end{align}

Since $\tilde u^+ = u$ holds on $K\cap \Omega^+$ we obtain
$$
|u-\tilde u^+|_{1,K} = |u-\tilde u^+|_{1,K\cap\, \Omega^-} \le |u^-|_{1,K\cap \, \Omega^-} +
|\tilde u^+|_{1,K\cap\, \Omega^-}.
$$
Applying H{\"o}lder's inequality with $p=\frac32$ and $q=3$, the imbedding of $H^1(K)$ into $L^6(K)$ (Note that
the imbedding constant can be bounded independently of $h$)
and \eqref{eq:measure-k-omega} one can bound $|u^-|_{1,K\cap \,\Omega^-}$ (and analogously
$|\tilde u^+|_{1,K\cap\, \Omega^-}$) by
\begin{align*}
|u^-|_{1,K\cap \,\Omega^-} &\le
|K\cap \Omega^-|^{\frac13} \|\nabla u^-\|_{0,6,K\cap\,\Omega^-}\\
&\le C\,h\,\|\nabla \tilde u^-\|_{0,6,K} \le C\,h\,|\tilde u^-|_{2,K}.
 \intertext{Hence}
 |u-\tilde u^+|_{1,K} &\le C\,h\, (|\tilde u^-|_{2,K} + |\tilde u^+|_{2,K}).
 \end{align*}
 Inserting this estimate into \eqref{eq:interp1} leads to
 \[
 |u-I_h u|_{1,K} \le C \frac{h^2}{\varrho_K}\, (|\tilde u^-|_{2,K} + |\tilde u^+|_{2,K}).
 \]
To prove the global interpolation error bound, we write
\begin{align*}
|u  - I_h u|_{1,\Omega}^2 &= \sum_{T\in \mathscr{T}_h\setminus \mathscr{T}_h^\gamma} |u  - I_h u|_{1,T}^2 + \sum_{T\in  \mathscr{T}_h^\gamma} \sum_{K\in \mathscr{T}_T^\gamma} |u  - I_h u|_{1,K}^2\\
&\le
Ch^2 \sum_{T\in \mathscr{T}_h\setminus \mathscr{T}_h^\gamma} |u|^2_{2,T} + C \sum_{T\in  \mathscr{T}_h^\gamma} \sum_{K\in \mathscr{T}_T^\gamma} \frac{h^2}{\varrho_K} (|\tilde u^-|^2_{2,K} + |\tilde u^+|^2_{2,K})\\
&\le
C\frac{h^2}{\varrho}\, (|\tilde u^-|^2_{2,\Omega} + |\tilde u^+|^2_{2,\Omega})\\
&\le C\frac{h^2}{\varrho}\, |u|^2_{2,\Omega^+\cup\,\Omega^-},
\intertext{where}
\varrho &=\min\{\varrho_K\colon K\in \mathscr{T}_T^\gamma, T\in \mathscr{T}_h^\gamma\}.\\[-2em]
\end{align*}
The calculation above indicates how the convergence rate can be improved in case $u\in W^{2,\infty}(\Omega^+\cup \Omega^-)$ observing that $|S_h^\gamma|\le Ch$ holds.
\end{proof}

\begin{rem}
It is classic in finite element theory to assume that the meshes are regular in the sense that
Condition \eqref{MeshRegularity} is satisfied for $\theta=0$. For the fitted meshes
$\mathscr T_h^\gamma$ one cannot guarantee that such a condition is satisfied.
To relax this constraint, we
assume here \eqref{MeshRegularity} for a $\theta\in [0,1)$ thus allowing a larger class of fitted
meshes than permitted by $\theta=0$.
\end{rem}

The following result gives the convergence rate for Problem \eqref{AFEM}.

\begin{thm}
\label{FirstErrorEstimate}
Assume that the family of fitted meshes $(\mathscr{T}_h^F)_h$ satisfies the regularity property
\eqref{MeshRegularity}. Then we have the error estimate
\begin{equation}
|u-u_h^F|_{1,\Omega} \le
\begin{cases}
Ch^{1-\theta}\,\|u\|_{2,\Omega^+\cup\,\Omega^-} &\text{if } u\in H^2(\Omega^+\cup\Omega^-),\\
C h\,\|u\|_{2,\infty,\Omega^+\cup\,\Omega^-}&\text{if } u\in W^{2,\infty}(\Omega^+\cup\Omega^-).
\end{cases}
\label{error-uhf}
\end{equation}
\end{thm}

\begin{proof}
We have from the triangle inequality
\begin{equation}
|u-u_h^F|_{1,\Omega} \le |u-\widehat u_h|_{1,\Omega} + |\widehat u_h-u_h^F|_{1,\Omega}.\label{TrIneq}
\end{equation}
To bound the first term on the right-hand side of \eqref{TrIneq}, we proceed as follows: Let us subtract
\eqref{Ah} from \eqref{Pb} and choose $v=u-\widehat u_h$. We have
$$
\int_\Omega (a\,\nabla u-a_h\nabla\widehat u_h)\cdot\nabla(u-\widehat u_h)\,dx = 0.
$$
Then
\begin{align*}
\int_\Omega &a_h |\nabla(u-\widehat u_h)|^2\,dx = -\int_\Omega (a-a_h)\,\nabla u\cdot
\nabla(u-\widehat u_h)\,dx\\
&=-\int_{\Omega\setminus S_h^\gamma}(a-a_h)\,\nabla u\cdot\nabla(u-\widehat u_h)\,dx -
\sum_{T\in \mathscr{T}_h^\gamma} \int_T (a-a_h)\,\nabla u\cdot\nabla(u-\widehat u_h)\,dx.
\end{align*}
The usual estimate for the interpolation error gives
\begin{align*}
\|a-a_h\|_{0,\infty,\Omega} &\le Ch\,(\|\tilde a\|_{1,\infty,\Omega^+_h}+\|\tilde a\|_{1,\infty,\Omega^-_h})\\
&\le Ch\,\|a\|_{1,\infty,\Omega^+\cup \Omega^-}.
\end{align*}
with a constant $C$ which only depends on a reference triangle, (see \cite{Ciarlet}, p.~124). Thus we obtain
\begin{equation}\label{eq:est-outside-band}
\bigg|\int_{\Omega\setminus S_h^\gamma}(a-a_h)\,\nabla u\cdot\nabla(u-\widehat u_h)\,dx\bigg| \le C\,h\,
\|a\|_{1,\infty,\Omega^+\cup \Omega^-}\, |u|_{1,\Omega\setminus S_h^\gamma}\,
|u-\widehat u_h|_{1,\Omega\setminus S_h^\gamma}.
\end{equation}
Next we consider a triangle $T\in \mathscr{T}_h^\gamma$ which we split  as
\[
T = (T\cap \Omega^+\cap\Omega_h^+)\cup (T\cap \Omega^-\cap\Omega_h^-)
\cup (T\cap \Omega^+\cap\Omega_h^-)\cup (T\cap \Omega^-\cap\Omega_h^+).
\]
As before, we obtain
\[
\bigg|\int_{T\cap \Omega^+\cap\Omega_h^+}(a-a_h)\,\nabla u\cdot \nabla(u-\widehat u_h)\,dx\bigg|
\le  Ch\, \|a\|_{1,\infty,\Omega^+\cup \Omega^-}|u|_{1,T\cap \Omega^+\cap\Omega_h^+}\,
|u-\widehat u_h|_{1,T\cap \Omega^+\cap\Omega_h^+}.
\]
Arguing as in the proof of Lemma~\ref{eq:interpolation-error}, the generalized H{\"o}lder inequality together with
\eqref{eq:measure-k-omega} yields the estimate
\begin{align*}
\bigg|\int_{T\cap \Omega^+\cap\Omega_h^-}& (a-a_h)\,\nabla u\cdot \nabla(u-\widehat u_h)\,dx\bigg|\\
&= \bigg|\int_{T\cap \Omega^+\cap\Omega_h^-} (a^+-a_h^-)\nabla u^+\cdot \nabla(u^+-\widehat u_h)\,dx\bigg|\\
&\le
C\, \|a\|_{0,\infty,\Omega}\, |T\cap \Omega^+\cap\Omega_h^-|^{1/3}\,\|\nabla u^+\|_{0,6,T\cap \Omega^+\cap\Omega_h^-}
\,\|\nabla(u^+-\widehat u_h)\|_{0,T\cap \Omega^+\cap\Omega_h^-}\\
&\le
C\, h\,\|a\|_{0,\infty,\Omega}\, |\tilde u^+|_{2,T}\,\|\nabla(u^+-\widehat u_h)\|_{0,T}.
\end{align*}
Analogous estimates hold with $+$ and $-$ interchanged. Collecting the four contributions to the
triangle $T$ one obtains
\begin{align*}
\bigg|\int_T&(a-a_h)\,\nabla u\cdot \nabla(u-\widehat u_h)\,dx\bigg|\\
 &\le C h\, (\|a\|_{0,\infty,\Omega} +\|a\|_{1,\infty,\Omega^+\cup\,\Omega^-})\\
 &\quad\times \big(|\tilde u^+|_{2,T}\,\|\nabla(u^+-\widehat u_h)\,
 \|_{0,T\cap\Omega^+} +  |\tilde u^-|_{2,T}\,\|\nabla(u^--\widehat u_h)\|_{0,T\cap\Omega^-}\big).
\end{align*}
Combining this estimate with \eqref{eq:est-outside-band} leads to
\begin{align*}
\int_\Omega &a_h\,|\nabla(u-\widehat u_h)|^2\,dx \le Ch\, \|a\|_{1,\infty,\Omega^+\cup\,\Omega^-}
|u|_{1,\Omega\setminus S_h^\gamma}|u- \widehat u_h|_{1,\Omega\setminus S_h^\gamma}\\
&+C\,h\, (\|a\|_{0,\infty,\Omega} + \|a\|_{1,\infty,\Omega^+\cup\,\Omega^-})\\
&\quad\times\sum_{T\in \mathscr{T}_h^\gamma}
\Big(|\tilde u^+|_{2,T}\|\nabla(u^+-\widehat u_h)\|_{0,T\cap\Omega^+}+
|\tilde u^-|_{2,T}\|\nabla(u^--\widehat u_h)\|_{0,T\cap\Omega^-}\Big)\\
&\le C\,h\,\|a\|_{1,\infty,\Omega^+\cup\,\Omega^-}\,|u|_{1,\Omega\setminus S_h^\gamma}\,
|u- \widehat u_h|_{1,\Omega\setminus S_h^\gamma} \\
&\qquad+C\,h\,(\|a\|_{0,\infty,\Omega} + \|a\|_{1,\infty,\Omega^+\cup\,\Omega^-})(|\tilde u^+|_{2,S_h^\gamma}
+ |\tilde u^-|_{2,S_h^\gamma})\,\|\nabla(u-\widehat u_h)\|_{0,S_h^\gamma}\\
&\le C\,h\,(\|a\|_{0,\infty,\Omega} + \|a\|_{1,\infty,\Omega^+\cup\,\Omega^-})\,
|u|_{2,\Omega^+\cup\,\Omega^-} \|\nabla(u-\widehat u_h)\|_{0,\Omega},
\end{align*}
which by \eqref{Prop-ah-3} implies
\begin{equation}\label{eq:est-u-hat-u}
|u-\widehat u_h|_{1,\Omega}\le C\,h\,(\|a\|_{0,\infty,\Omega} + \|a\|_{1,\infty,\Omega^+\cup\,\Omega^-})\,|u|_{2,\Omega^+\cup\,\Omega^-}.
\end{equation}
To bound the norm $|\widehat u_h-u_h^F|_{1,\Omega}$, we have from problems \eqref{Ah} and \eqref{AFEM},
$$
\int_\Omega a_h\nabla(\widehat u_h-u_h^F)\cdot\nabla v\,dx = 0\qquad\forall\ v\in W_h.
$$
Standard finite element approximation theory combined with \eqref{Prop-ah-1}--\eqref{Prop-ah-2} gives
\begin{equation}
|\widehat u_h-u_h^F|_{1,\Omega} \le C\,\inf_{v\in W_h}|\widehat u_h-v|_{1,\Omega},
\label{AbstractBound}
\end{equation}
which together with \eqref{eq:est-u-hat-u} implies
\begin{align*}
|\widehat u_h-u_h^F|_{1,\Omega} &\le C\,|\widehat u_h - I_h u|_{1,\Omega}\\
&\le C\,|\widehat u_h - u|_{1,\Omega} + C\,|u - I_h u|_{1,\Omega}\\
&\le
C\,h\,|u|_{2,\Omega^+\cup\,\Omega^-} + C\,|u - I_h u|_{1,\Omega}.
\end{align*}
The interpolation error is bounded using \eqref{InterpolationErrorGlobal-1} or \eqref{InterpolationErrorGlobal-2}.
\end{proof}

\section{A hybrid approximation}
The method presented in the previous section has proven its efficiency as numerical tests will show
in the last section. In  more elaborate problems like time dependent
or nonlinear problems where the interface $\gamma$ is a moving front, the subtriangulation $\mathscr T_h^\gamma$
moves within iterations and then the matrix structure has to be frequently modified. To remedy to this difficulty, we resort to
a hybridization of the added unknowns. More specifically, the added discrete space $X_h$ is replaced by a nonconforming
approximation space. In addition, a Lagrange multiplier is used to compensate this inconsistency.
The hybridization enables to locally eliminate the added unknowns in each triangle
$T\in\mathscr T_h^\gamma$. In the sequel we fix an orientation for the interface $\gamma$.
This induces an orientation of the normals to the edges $e\in \mathscr{E}_h^\gamma$ by
following the interface in the positive direction. The jump of a function $v$ across an
edge $e\in \mathscr{E}_h^\gamma$ can then be defined as
$$
[v]_e(x) := \lim_{s\to 0, s>0}v(x+s n(x)) - \lim_{s\to 0, s<0}v(x+s n(x))\equiv v^+(x)-v^-(x),\qquad x\in e,
$$
where $n$ is the  unit normal to $e$.

To develop this method, we start by defining an ad-hoc formulation for the solution $\widehat u_h$ of \eqref{Ah}.
Let us define the spaces
\begin{align*}
\widehat Z_h &:= H^1_0(\Omega) + \widehat Y_h,\\
\widehat Y_h &:=\{v\in L^2(\Omega);\ v_{|\Omega\setminus S_h^\gamma}=0,
\ v_{|T}\in H^1(T)\ \forall\ T\in\mathscr T_h^\gamma,\\
&\qquad [v]=0\text{ on }e,\ \forall\ e\in\mathscr E_h\setminus\mathscr E_h^\gamma\},\\
\widehat Q_h &:= \prod_{e\in\mathscr E_h^\gamma}H^{-\frac 12}_{00}(e),
\end{align*}
where $H^{-\frac 12}_{00}(e)$ is the dual space of the trace space
$$
H^{\frac 12}_{00}(e) := \{v_{|e};\ v\in H^1(T),\ e\in\mathscr E_T,\ v=0\text{ on }d\quad \forall\ d\in\mathscr E_T, d\ne e\}.
$$
We remark that the jumps $[v]$ for $v\in\widehat Z_h$ can be interpreted in $H^{\frac 12}_{00}(e)$ for $e\in\mathscr E_h^\gamma$. This is due to the fact that $v\in H^1(T)$ for all $T\in\mathscr T_h$, that for every $e\in\mathscr E_h^\gamma$, the jump of $v$ lies in $H^{\frac 12}(e)$ and vanishes at the endpoints of $e$ as well as on at least two adjacent edges. This motivates the choice of $\widehat Q_h$.

The elements of  $\widehat Q_h$  will be referred to by $\mu = (\mu_e)_{e\in \mathscr E_h^\gamma}$.
We endow $\widehat Z_h$ with the broken norm
$$
\|u\|_{\widehat Z_h} = (\sum_{T\in \mathscr T_h} |u|^2_{1,T})^{1/2}.
$$
On $\widehat Q_h$ we use the norm
$$
\|\mu\|_{\widehat Q_h} = \Big(\sum_{e\in \mathscr E_h^\gamma} \|\mu_e\|^2_{H^{-\frac 12}_{00}(e)}\Big)^{\frac 12}
:= \Bigg(\sum_{e\in \mathscr E_h^\gamma}
\bigg(\sup_{v\in H^{\frac 12}_{00}(e)\setminus\{0\}} \frac{\int_e \mu_e v\,ds}{\|v\|_{H^{\frac 12}_{00}(e)}}\bigg)^2\Bigg)^{\frac 12}.
$$
Above, the integrals over edges $e$ are to be interpreted as  duality pairings between $H^{-\frac 12}_{00}(e)$ and $H^{\frac 12}_{00}(e)$.
We mention that the broken norm in $\widehat Z_h$ reflects the fact that $\widehat Z_h$ is not a subspace of $H_0^1(\Omega)$.

Next we  define the variational problem,
\begin{alignat}{2}
&\text{Find }(\widehat u_h^H,\widehat\lambda_h)\in \widehat Z_h\times \widehat Q_h\text{ such that:}\nonumber\\
&\sum_{T\in\mathscr T_h}\int_T a_h\,\nabla \widehat u_h^H\cdot\nabla v\,dx - \sum_{e\in\mathscr E_h^\gamma}\int_e\widehat\lambda_h\, [v]\,ds = \int_\Omega fv\,dx
&&\qquad\forall\ v\in\widehat Z_h,\label{PbLambda-1}\\
&\sum_{e\in\mathscr E_h^\gamma}\int_e\mu\, [\widehat u_h^H]\,ds = 0&&\qquad\forall\ \mu\in \widehat Q_h.\label{PbLambda-2}
\end{alignat}
The saddle point problem \eqref{PbLambda-1}--\eqref{PbLambda-2} indicates that the continuity of $\widehat u_h$ across the edges of $\mathscr E_h^\gamma$ is enforced by a Lagrange multiplier technique.

\begin{thm}
\label{Existence-uH}
Problem \eqref{PbLambda-1}--\eqref{PbLambda-2} has a unique solution
$(\widehat u_h^H,\widehat\lambda_h)\in \widehat Z_h\times \widehat Q_h$. Moreover, we have
$\widehat u_h^H=\widehat u_h$ and the following estimate holds
\begin{equation}
\|\widehat u_h^H\|_{\widehat Z_h} + \|\widehat\lambda_h\|_{\widehat Q_h} \le C\,\|f\|_{0,\Omega},\label{estim-ul}
\end{equation}
with a constant $C$ which is independent of $h$.
\end{thm}

\begin{proof}
Problem \eqref{PbLambda-1}--\eqref{PbLambda-2} can be put in the standard variational form
$$
\left\{
\begin{aligned}{}
&\mathscr A(\widehat u^H_h,v) + \mathscr B(v,\widehat\lambda_h) = (f,v)&&\qquad\forall\ v\in \widehat Z_h,\\
&\mathscr B(\widehat u^H_h,\mu) = 0&&\qquad\forall\ \mu\in \widehat Q_h,
\end{aligned}
\right.
$$
where
\begin{align*}
&\mathscr A(u,v) = \sum_{T\in\mathscr T_h}\int_T a_h\,\nabla u\cdot\nabla v\,dx,\\
&\mathscr B(v,\mu) = - \sum_{e\in\mathscr E_h^\gamma}\int_e \mu\,[v]\,ds,\\
&(f,v) = \int_\Omega fv\,dx.
\end{align*}
The bilinear form $\mathscr A$ is clearly continuous and coercive on the space
$\widehat Z_h\times\widehat Z_h$.
The bilinear form $\mathscr B$ is also continuous on $\widehat Z_h\times\widehat Q_h$.

Next we verify that $\mathscr B$ satisfies the inf-sup condition, i.e. there exists $\delta >0$
such that for every $\lambda\in \widehat Q_h$ there exists $v_\mu\in \widehat Z_h$ such that
$$
\mathscr B(v_\mu,\mu) \ge \delta\, \|v_\mu\|_{\widehat Z_h}\|\mu\|_{\widehat Q_h}
$$
{\em i.e.}
\begin{equation}\label{eq:inf-sup-semidiscr}
\sum_{e\in \mathscr{E}_h^\gamma} \int_e \mu_e[v_\mu]\,ds \ge \delta\,\|v_\mu\|_{\widehat Z_h}\|\mu\|_{\widehat Q_h}
\end{equation}
holds.

Given $\mu =(\mu_e)_{e\in\mathscr E^\gamma_h}\in \widehat Q_h$ and an edge $e\in\mathscr{E}_h^\gamma$
choose a triangle $T\in \mathscr{T}_h^\gamma$ which has $e$ as one of its edges.
Define $v_T\in H^1(T)$ as the solution of
\begin{equation}
\left\{
\begin{aligned}{}
&\Delta v = 0&&\qquad\text{in } T,\\
&\frac{\partial v}{\partial n} = \mu_e &&\qquad\text{on } e,\\
&v = 0 &&\qquad\text{on } \partial T\setminus e,
\end{aligned}
\right.
\label{Pb-mu-v}
\end{equation}
which is equivalent to
$$
\int_T \nabla v\cdot\nabla \varphi\,dx = \int_e\mu_e \varphi\,ds\quad\text{for }\varphi\in H^1_e(T)
$$
where
$$
H^1_e(T) = \{\varphi\in H^1(T);\ \varphi=0 \text{ on }\partial T\setminus e \}.
$$
By Green's theorem we obtain
\begin{align*}
\|\mu_e\|_{-1/2,e} &= \Big\|\frac{\partial v_T}{\partial n}\Big\|_{-1/2,e}\le \|\nabla v_T\|_{0,T},\\
\int_e\mu_e v_T\, ds &= \int_T |\nabla v_T|^2\,dx,
\end{align*}
which implies
\[
\|\mu_e\|_{-1/2,e}^2 \le \int_T |\nabla v_T|^2\,dx = \int_e\mu_e v_T\, ds.
\]
Let $\chi_T$ denote the characteristic function of $T$ and define
\[
v_\mu = \sum_{T\in \mathscr{T}_h^\Gamma} \chi_Tv_T.
\]
Since there are as many edges in $\mathscr{E}_h^\gamma$ as triangles in $\mathscr{T}_h^\gamma$ then
$[v_\mu] = v_T$ holds for every edge $e\in \mathscr{E}_h^\gamma$. Hence we obtain
\[
\|\mu\|^2_{\widehat Q_h} = \sum_{e\in \mathscr{E}_h^\gamma} \|\mu_e\|_{-1/2,e}^2 \le
\sum_{T\in\mathscr{T}_h^\gamma} \|\nabla v_T\|^2_{0,T} = \sum_{e\in \mathscr{E}_h^\gamma}\int_e
\mu_e [v_\mu]\,ds.
\]
Furthermore,
\[
\|v_\mu\|^2_{\widehat Z_h} = \sum_{T\in\mathscr{T}_h}\|\nabla v_T\|^2_{0,T} =
\sum_{T\in\mathscr{T}_h^\gamma}\|\nabla v_T\|^2_{0,T}
\]
holds. This implies
\[
\|\mu\|^2_{\widehat Q_h}\|v_\mu\|^2_{\widehat Z_h} \le \Big(\sum_{T\in\mathscr{T}_h^\gamma}
\|\nabla v_T\|^2_{0,T}\Big)^2 = \mathscr B(v_\mu,\mu)^2.
\]
Adjusting the sign of $v_\mu$ this is equivalent to \eqref{eq:inf-sup-semidiscr} with $\delta = 1$. The estimate \eqref{estim-ul} is a direct
consequence of \eqref{eq:inf-sup-semidiscr}.

Now, it is clear from \eqref{PbLambda-2} that
$$
[\widehat u_h^H]=0\quad\text{on e},\quad\forall\ e\in\mathscr T_h^\gamma.
$$
This implies that $\widehat u_h^H\in H^1_0(\Omega)$. Choosing a test function $v\in H^1_0(\Omega)$ in \eqref{PbLambda-1},
we find that $\widehat u_h^H$ is a solution to Problem \eqref{Ah}, and then $\widehat u^H_h=\widehat u_h$.
The interpretation of $\widehat\lambda_h$ is simply obtained by the Green's formula.
\end{proof}

We are now able to present a numerical method to solve the interface problem. This one is simply derived
as a finite element method to solve the saddle point problem \eqref{PbLambda-1}--\eqref{PbLambda-2}. We consider for this end
a piecewise constant approximation of the Lagrange multiplier. Let us define the finite dimensional spaces,
\begin{align*}
Z_h &:= V_h + Y_h,\\
Y_h &:=\{v\in L^2(\Omega);\ v_{|\Omega\setminus S_h^\gamma}=0,\ v_{|K}\in
P_1(K)\ \forall\ K\in\mathscr T_T^\gamma,\ \forall\ T\in\mathscr T_h^\gamma,\\
&\qquad [v]=0\text{ on }e,\ \forall\ e\in\mathscr E_h\setminus\mathscr E_h^\gamma\},\\
Q_h &:= \big\{\mu\in \prod_{e\in\mathscr E_h^\gamma}L^2(e);\ \mu_{|e}=\text{const.}\quad \forall\ e\in\mathscr E_h^\gamma\big\}.
\end{align*}

The hybrid finite element approximation is given by the following problem:
\begin{alignat}{2}
&\text{Find }(u^H_h,\lambda_h)\in Z_h\times Q_h\text{ such that:}\nonumber\\
&\sum_{T\in\mathscr T_h}\int_T a_h\,\nabla u^H_h\cdot\nabla v\,dx - \sum_{e\in\mathscr E_h^\gamma}\int_e\lambda_h\, [v]\,ds = \int_\Omega fv\,dx
&&\qquad\forall\ v\in Z_h,\label{PbLambda-h-1}\\
&\sum_{e\in\mathscr E_h^\gamma}\int_e\mu\, [u^H_h]\,ds = 0&&\qquad\forall\ \mu\in Q_h.\label{PbLambda-h-2}
\end{alignat}
Let us give some additional remarks before proving convergence properties of this method.

1. The matrix formulation of the method has the following form
\begin{equation}
\begin{pmatrix}
A & C & 0 \\ C^T & D & B \\ 0 & B^T &  0
\end{pmatrix}
\begin{pmatrix}
\widetilde u \\ \widetilde v \\ \widetilde\lambda
\end{pmatrix}
=
\begin{pmatrix}
b \\ c \\ 0
\end{pmatrix},
\label{LS}
\end{equation}
where the vector $\widetilde u$ contains the values of $u_h^H$ at nodes of the mesh $\mathscr T_h$, {\em i.e.} components
of $u_h^H$ in the Lagrange basis of $V_h$, $\widetilde v$ contains the components of $u_h^F$ in the basis of
$Y_h$, and $\widetilde \lambda$ has as components the values of $\lambda_h$  on the edges of $\mathscr E_h^\gamma$. There is clearly no simple method to
eliminate off diagonal blocks  in the system \eqref{LS} in order to decouple the variables. More specifically,
our aim is to eliminate the unknowns $\widetilde v$.

2. The method must be viewed in the context of an iterative process like the Uzawa method, where
the Lagrange multiplier $\lambda_h$ is decoupled from the primal variable $u_h^H$. In such situations, each iteration step
consists in solving an elliptic problem with a given $\lambda_h$. Let us recall that, due to the local feature
of the basis functions of nodes on edges of $\mathscr E_h^\gamma$, the unknowns associated to these nodes
can be eliminated at the element level. This is a basic issue in our method.

3. We point out that equation \eqref{PbLambda-h-2} entails
\begin{equation}
[u_h^H]=0\quad\text{on }e,\quad\forall\ e\in\mathscr T_h^\gamma.\label{Cont-uhH}
\end{equation}
This follows from the fact  that  $u_h^H$ is an affine function on each edge of $\mathscr T_h^\gamma$. This implies that actually $u_h^H\in W_h$. Choosing $v\in W_h$ in \eqref{PbLambda-1}  we find
$$
\int_\Omega a_h\nabla u_h^H\cdot\nabla v\,dx = \int_\Omega fv\,dx.
$$
This yields $u_h^H= u_h^F$.

\section{Convergence analysis}

This section is devoted to the proof of existence, uniqueness and stability of the solution
of \eqref{PbLambda-h-1}--\eqref{PbLambda-h-2} as well as its convergence to Problem \eqref{PbLambda-1}--\eqref{PbLambda-2}.

For this result we need a localized quasi-uniformity of the mesh. More precisely, we assume that
\begin{equation}
|e| \ge Ch\qquad\forall\ e\in\mathscr E_h^\gamma.\label{QU}
\end{equation}
In addition, we make the following assumption:
\begin{equation}
\begin{aligned}{}
&\text{The distance of the intersection point of $\gamma$ with any edge $e\in\mathscr E_h^\gamma$}\\
&\text{to the endpoints of $e$ can be bounded from below by $\delta h$, where $\delta$ is}\\
&\text{independent of $h$.}
\label{Hyp}
\end{aligned}
\end{equation}
Although this assumption appears to be quite restrictive, numerical tests have shown that it can be actually ignored in applications.

%

\begin{thm}
\label{EU:Discrete}
Assume that the family of meshes $(\mathscr T_h)_h$ satisfies Property \eqref{QU}.
Then Problem \eqref{PbLambda-h-1}--\eqref{PbLambda-h-2} has a unique solution. Moreover, we have the bound
\begin{equation}
\|u^H_h\|_{\widehat Z_h} + \|\lambda_h\|_{Q_h} \le C\,\|f\|_{0,\Omega},
\label{Stability}
\end{equation}
where the constant $C$ is independent of $h$.
\end{thm}

\begin{proof}
It is clearly sufficient to prove the inf-sup condition (see for instance Brezzi-Fortin \cite{BF}):
\begin{equation}
\sup_{v_h\in Z_h\setminus\{0\}}\dfrac{\sum_{e\in\mathscr E_h^\gamma}\int_e\mu_h\,[v_h]\,ds}
{\|v_h\|_{\widehat Z_h}\,\|\mu_h\|_{\widehat Q_h}}\ge \beta > 0
\qquad\forall\ \mu_h\in Q_h.
\label{LBB}
\end{equation}
In the following, for each triangle $T\in\mathscr T_h^\gamma$, we shall denote
by $e^+_T$ ({\sl resp.} $e^-_T$) the edge where $\gamma$ enters $T$ ({\sl resp.} leaves $T$), and by $\tilde e_T$
the remaining edge of $T$ (see Figure \eqref{Fig-2}). Recall that we fixed an orientation for $\gamma$.
\begin{figure}[!ht]
\centering
\includegraphics*[height=5cm,viewport=100 520  500 720]{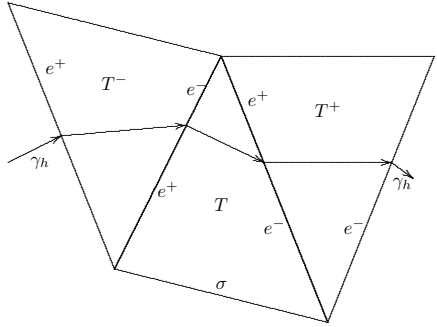}
\caption{\label{Fig-2} Definition of $e^+_T$, $e^-_T$, and $\tilde e_T$.}
\end{figure}

Let $\mu_h\in Q_h$, and let $v\in\widehat Z_h$ be the function given by Problem \eqref{Pb-mu-v}.

We define a function $v_h\in Z_h$ by
\begin{equation}
\left\{
\begin{aligned}{}
&v_{h|T} = 0 &&\qquad\forall\ T\in\mathscr T_h\setminus\mathscr T_h^\gamma,\\
&\int_e v_h\,ds = \int_e v\,ds &&\qquad\forall\ e\in\mathscr E_T,\ \forall\ T\in\mathscr T_h^\gamma.
\end{aligned}
\right.
\label{def-vh}
\end{equation}
The gradient of $v_h$ can be expressed in $T\in\mathscr T_h^\gamma$ by
$$
\nabla v_{h|T} = \frac 2{|e^-_T|}\bigg(\int_{e^-_T}v\,ds\bigg)\,\nabla\varphi_{e^-_T}
+ \frac 2{|e^+_T|}\bigg(\int_{e^+_T}v\,ds\bigg)\,\nabla\varphi_{e^+_T},
$$
where $\varphi_{e^+_T}$ ({\sl resp.} $\varphi_{e^-_T}$) is the basis function of $Z_h$ associated to the
added node on $e^+_T$ ({\sl resp.} $e^-_T$).
Then by using \eqref{QU} and the Cauchy-Schwarz inequality, we get for each $T\in\mathscr T_h^\gamma$,
\begin{align}
\|\nabla v_h\|_{0,T} &= C_1\,h^{-1}\bigg|\int_{e^-_T}v\,ds\bigg|\,\|\nabla\varphi_{e^-_T}\|_{0,T}
+ C_2\,h^{-1}\bigg|\int_{e^+_T}v\,ds\bigg|\,\|\nabla\varphi_{e^+_T}\|_{0,T}\nonumber\\
&\le C_3\,h^{-\frac 12}\,\big(\|v\|_{0,e_T^-}\|\nabla\varphi_{e^-_T}\|_{0,T}
+ \|v\|_{0,e_T^+}\,\|\nabla\varphi_{e^+_T}\|_{0,T}\big).\label{ident-1}
\end{align}
The trace inequality (see \cite{Arnold}, eq. (2.5)) and the Poincar\'e inequality
owing to $v=0$ on $\tilde e_T$, yield for $T\in\mathscr T_h^\gamma$,
\begin{equation}
\|v\|_{0,e^\pm_T}
\le C_4\,\big(h^{-\frac 12}\|v\|_{0,T}+h^{\frac 12}\|\nabla v\|_{0,T}\big)
\le C_5\,h^{\frac 12}\|\nabla v\|_{0,T}.\label{ident-2}
\end{equation}
On the other hand, Assumption \eqref{Hyp} implies the uniform boundedness of
$\|\nabla\varphi_{e^\pm_T}\|_{0,T}$. From \eqref{ident-1} and \eqref{ident-2} we obtain then
$$
\|\nabla v_h\|_{0,T} \le C_6\,\|\nabla v\|_{0,T}.
$$
Using the inf-sup condition \eqref{eq:inf-sup-semidiscr} and \eqref{def-vh}, we finally obtain
\begin{align*}
\|\mu_h\|_{\widehat Q_h}\|v_h\|_{\widehat Z_h}
&\le C_6\,\|\mu_h\|_{\widehat Q_h}\|v\|_{\widehat Z_h}\\
&\le C_7 \sum_{e\in\mathscr E_h^\gamma}\int_e \mu_h\,[v]\,ds\\
&= C_7\,\sum_{e\in\mathscr E_h^\gamma}\int_e \mu_h\,[v_h]\,ds.
\end{align*}
Finally, obtaining the estimate \eqref{Stability} is a classical task that we skip here.
\end{proof}

We now prove the main convergence result.

\begin{thm}
\label{Convergence-u}
Assume hypotheses \eqref{MeshRegularity} and \eqref{QU} are satisfied, then
there exists a constant $C$, independent of $h$, such that
$$
\|u-u_h^H\|_{\widehat Z_h} \le
\begin{cases}
Ch^{1-\theta}\,|u|_{2,\Omega^+\cup\Omega^-}&\text{if }u\in H^2(\Omega^+\cup\Omega^-),\\
Ch\,\|u\|_{2,\infty,\Omega^+\cup\Omega^-}&\text{if }u\in W^{2,\infty}(\Omega^+\cup\Omega^-).
\end{cases}
$$
\end{thm}

\begin{proof}
From classical theory of saddle point problems (see \cite{GiraultRaviart}, p. 114),
we obtain from Theorem \ref{EU:Discrete},
\begin{equation}
\|\widehat u_h^H-u_h^H\|_{\widehat Z_h} + \|\widehat\lambda_h-\lambda_h\|_{\widehat Q_h}\le C\,
\Big(\inf_{v\in Z_h}\|\widehat u_h^H-v\|_{\widehat Z_h}+\inf_{\mu\in Q_h}\|\widehat\lambda_h-\mu\|_{\widehat Q_h}\Big).
\label{AbstractErrorEstimate-1}
\end{equation}
Furthermore, using (Braess \cite{Braess}, Theorem 4.8), Property \eqref{Cont-uhH}
implies that Estimate \eqref{AbstractErrorEstimate-1} can be improved, for the error on $\widehat u_h^H$ by
\begin{equation}
\|\widehat u_h^H-u_h^H\|_{\widehat Z_h}\le C\,\inf_{v\in Z_h}\|\widehat u_h^H-v\|_{\widehat Z_h}.
\label{AbstractErrorEstimate-2}
\end{equation}
To bound the right-hand side, we choose $v=I_hu$, where $I_h$ is the previously defined
Lagrange interpolant in $Z_h$. Since $\widehat u_h^H=\widehat u_h$ (see Theorem \ref{Existence-uH}),
then by using \eqref{InterpolationErrorGlobal-1} and \eqref{eq:est-u-hat-u},
\begin{align*}
\|\widehat u_h^H-I_hu\|_{\widehat Z_h} &= \|\widehat u_h-I_hu\|_{\widehat Z_h}\\
&\le \|u-I_h u\|_{\widehat Z_h} + \|u-\widehat u_h\|_{\widehat Z_h}\\
&\le C_1\,h^{1-\theta}\,|u|_{2,\Omega^+\cup\Omega^-} + C_2\,h\,|u|_{2,\Omega^+\cup\Omega^-}.
\end{align*}
If $u\in W^{2,\infty}(\Omega^+\cup\Omega^-)$, then Estimate \eqref{InterpolationErrorGlobal-2} yields
\begin{equation}
\|\widehat u_h^H-I_hu\|_{\widehat Z_h} \le C\,h\,\|u\|_{2,\infty,\Omega^+\cup\Omega^-}.
\label{ErrorInterp}
\end{equation}
\end{proof}

\begin{rem}
As it was previously mentioned, we know that if Problem \eqref{PbLambda-1}--\eqref{PbLambda-2} has a unique
solution $(u_h^H,\lambda_h)$ then $u_h^H=u_h^F$ where $u_h^F$ is the solution of Problem
\eqref{AFEM} and therefore the error estimate \eqref{error-uhf} holds. Consequently, Theorem \ref{Convergence-u}
can be simply proven by obtaining a nonuniform inf-sup condition ({\em i.e.} \eqref{LBB} with $\beta=\beta(h))$.
This can be achieved without assuming \eqref{Hyp}. In this case, no error estimate is to be expected for
the Lagrange multiplier.
\end{rem}

Finally, since the Lagrange multiplier $\widehat\lambda_h$ can be interpreted in terms of $\widehat u_h$
(see Theorem \ref{Existence-uH}), it is interesting to see how good is its approximation $\lambda_h$.
Let, for this, $E_h$ denote the set
$$
E_h := \prod_{e\in\mathscr E_h^\gamma}e.
$$

\begin{thm}
Under the same hypotheses as in Theorem \ref{Convergence-u}, we have the following error bounds
$$
\|\widehat\lambda_h-\lambda_h\|_{\widehat Q_h}\le
\begin{cases}
C(h^{1-\theta} + h^{\frac 12})\,\Big(|u|_{2,\Omega^+\cup\Omega^-}+\|\lambda\|_{0,E_h}\Big)
&\text{if }\widehat\lambda_h\in L^2(E_h),\\
Ch^{1-\theta}\,\Big(|u|_{2,\Omega^+\cup\Omega^-}+\|\lambda\|_{\frac 12,E_h}\Big)
&\text{if }\widehat\lambda_h\in H^{\frac 12}(E_h).
\end{cases}
$$
\end{thm}

\begin{proof}
We use the abstract error bound \eqref{AbstractErrorEstimate-1}. Let, for $e\in\mathscr E_h^\gamma$,
$$
\lambda_e := \frac 1{|e|}\int_e\widehat\lambda_h\,ds.
$$
Using Lemma 7 in Girault--Glowinski \cite{GG}, we obtain the bound
$$
\|\widehat\lambda_h-\lambda_e\|_{H^{-\frac 12}_{00}(e)}\le Ch^{\frac 12}\|\widehat\lambda_h\|_{0,e}
\qquad\text{if }\widehat\lambda_h\in L^2(e),
$$
and
$$
\|\widehat\lambda_h-\lambda_e\|_{H^{-\frac 12}_{00}(e)}\le Ch\,\|\widehat\lambda_h\|_{\frac 12,e}
\qquad\text{if }\widehat\lambda_h\in H^{\frac 12}(e),
$$
Combining these bounds with \eqref{AbstractErrorEstimate-1}, \eqref{AbstractErrorEstimate-2}
and \eqref{ErrorInterp} achieves the proof.
\end{proof}

\section{A numerical test}
To test the efficiency and accuracy of our method, we present in this section a numerical test.
We consider an exact radial solution and test convergence rates in various norms.

Let $\Omega$ denote the square $\Omega=(-1,1)^2$ and let the function $a$ be given by
$$
a(x) = \begin{cases}
\alpha &\text{if }|x|<R_1,\\
\beta &\text{if }|x|\ge R_1,
\end{cases}
$$
where $\alpha,\beta>0$. We test the exact solution
$$
u(x) =
\begin{cases}
\dfrac 1{4\alpha}\,(R_1^2-|x|^2) + \dfrac 1{4\beta}(R^2_2-R^2_1) &\text{if }|x|<R_1,\\
\dfrac 1{4\beta}\,(R^2_2-|x|^2) &\text{if }|x|\ge R_1.
\end{cases}
$$
We choose $R_1=0.5$ and $R_2=\sqrt{2}$.
The function $f$ and Dirichlet boundary conditions are determined according to this choice. Note that unlike
the presented model problem, we deal here with non homogeneous boundary conditions but this cannot
affect the obtained results.

The finite element mesh is made of $2N^2$ equal triangles. According to the definition of $a$,
the interface $\gamma$ is given by the circle of center $0$ and radius $R_1$. The error is measured
in the following discrete norms:
\begin{align*}
&\|e\|_{0,h} := \bigg(\frac 1M\,\sum_{i=1}^M (u(x_i)-u_h(x_i))^2\bigg)^{\frac 12},\\
&\|e\|_{0,\infty} := \max_{1\le i\le M}|u(x_i)-u_h(x_i)|,\\
&\|e\|_{1,h} := \bigg(\sum_{T\in\mathscr T_h} \int_T|I_h(\nabla u)(x)-\nabla u_h|^2\bigg)^{\frac 12},
\end{align*}
where $x_i$ are the mesh nodes, $M$ is the total number of nodes, and $I_h$ is the piecewise linear interpolant.
We denote in the sequel by $p$ the ratio
$\alpha/\beta$.
Table 1 presents convergence rates for the standard $P_1$ finite element method using the unfitted mesh
\eqref{SFEM} with the choice $p=1/10$.
\bigskip
\begin{center}
\begin{tabular}{|c|c|c|c|c|c|c|}
\hline
 $h^{-1}$ & $\|e\|_{0,h}$ & Rate & $\|e\|_{0,\infty}$ & Rate & $\|e\|_{1,h}$ & Rate \\
\hline
$\phantom{1}10$ & $1.40\times 10^{-2}$ &      & $2.02\times 10^{-2}$ &      & $6.28\times 10^{-2}$ & \\
\hline
$\phantom{1}20$ & $6.78\times 10^{-3}$ & $1.05$ & $1.09\times 10^{-2}$ & $0.89$ & $5.23\times 10^{-2}$ & $0.26$ \\
\hline
$\phantom{1}40$ & $3.61\times 10^{-3}$ & $0.91$ & $5.81\times 10^{-3}$ & $0.91$ & $3.68\times 10^{-2}$ & $0.51$ \\
\hline
$\phantom{1}80$ & $1.83\times 10^{-3}$ & $0.98$ & $3.06\times 10^{-3}$ & $0.92$ & $2.56\times 10^{-2}$ & $0.52$ \\
\hline
$160$ & $9.44\times 10^{-4}$ & $0.95$ & $1.55\times 10^{-3}$ & $0.98$ & $1.82\times 10^{-2}$ & $0.49$ \\
\hline
\end{tabular}\\
\medskip
Table 1. Convergence rates for a standard (unfitted) finite element method.
\end{center}

\medskip
As expected, numerical experiments show poor convergence behavior. Let us consider now the results
obtained by the present method, \emph{i.e.} \eqref{PbLambda-1}--\eqref{PbLambda-2} or equivalently
\eqref{AFEM}. We obtain for $p=1/10$ and $p=1/100$ the convergence rates illustrated in Tables 1 and 2
respectively.

\bigskip
\begin{center}
\begin{tabular}{|c|c|c|c|c|c|c|}
\hline
 $h^{-1}$ & $\|e\|_{0,h}$ & Rate & $\|e\|_{0,\infty}$ & Rate & $\|e\|_{1,h}$ & Rate \\
\hline
$\phantom{1}10$ & $3.45\times 10^{-3}$ &        & $4.25\times 10^{-3}$ &      & $1.75\times 10^{-2}$ & \\
\hline
$\phantom{1}20$ & $8.18\times 10^{-4}$ & $2.1$  & $1.72\times 10^{-3}$ & $1.3$ & $6.87\times 10^{-3}$ & $1.3$ \\
\hline
$\phantom{1}40$ & $1.70\times 10^{-4}$ & $2.3$ & $5.22\times 10^{-4}$ & $1.7$ & $2.81\times 10^{-3}$ & $1.3$ \\
\hline
$\phantom{1}80$ & $3.94\times 10^{-5}$ & $2.1$ & $1.64\times 10^{-4}$ & $1.7$ & $1.02\times 10^{-3}$ & $1.5$ \\
\hline
$160$ & $8.57\times 10^{-6}$ & $2.2$ & $4.89\times 10^{-5}$ & $1.7$ & $3.59\times 10^{-4}$ & $1.5$ \\
\hline
\end{tabular}\\
\medskip
Table 2. Convergence rates for the hybrid finite element method with $p=1/10$.
\end{center}
\bigskip
\begin{center}
\begin{tabular}{|c|c|c|c|c|c|c|}
\hline
 $h^{-1}$ & $\|e\|_{0,h}$ & Rate & $\|e\|_{0,\infty}$ & Rate & $\|e\|_{1,h}$ & Rate \\
\hline
$\phantom{1}10$ & $3.26\times 10^{-3}$ &      & $4.07\times 10^{-3}$ &      & $1.69\times 10^{-2}$ & \\
\hline
$\phantom{1}20$ & $7.91\times 10^{-4}$ & $2.0$ & $1.74\times 10^{-3}$ & $1.2$ & $6.65\times 10^{-3}$ & $1.3$ \\
\hline
$\phantom{1}40$ & $1.72\times 10^{-4}$ & $2.2$ & $5.47\times 10^{-4}$ & $1.7$ & $2.72\times 10^{-3}$ & $1.3$ \\
\hline
$\phantom{1}80$ & $4.01\times 10^{-5}$ & $2.1$ & $1.74\times 10^{-4}$ & $1.6$ & $9.88\times 10^{-4}$ & $1.5$ \\
\hline
$160$ & $8.82\times 10^{-6}$ & $2.2$ & $5.22\times 10^{-5}$ & $1.7$ & $3.50\times 10^{-4}$ & $1.5$ \\
\hline
\end{tabular}\\
\medskip
Table 3. Convergence rates for the hybrid finite element method with $p=1/100$.
\end{center}

\medskip
Tables 2 and 3 show convergence rates that are even better than the theoretical results.
This is probably due to the choice of a discrete norm but may also be due to a superconvergence
phenomenon.
Rates for the $L^2$--norm give also good behavior. For the $L^\infty$--convergence rate,
we can note that we dot retrieve the second order obtained for a continuous coefficient problem.
However, these rates are better ($1.5$ rather than $1$) than the ones obtained for a standard finite element
method and, moreover, the error values are significantly lower in our case. It is in addition
remarkable that the error values depend very weekly on $p$ but the convergence rates are
independent of this value.

\section{Concluding remarks}

We have presented an optimal rate finite element method to solve interface problems with unfitted meshes.
The main advantage of the method is that the added unknowns that deal with the interface singularity
do not modify the matrix structure. This feature enables using the method in more complex situations
like in problems with moving interfaces. The price to pay for this is the use of a Lagrange multiplier
that adds an unknown on each edge that cuts the interface. This drawback can be easily
removed by using an iterative method such as the classical Uzawa method or more elaborate methods
like the Conjugate Gradient. The good properties of the obtained saddle point problem enable
choosing among a wide variety of dedicated methods. This topic will be addressed in a future work.
Let us also mention that the present finite element method does not specifically address problems
with large jumps in the coefficients. These ones are in addition ill conditioned and this drawback
is not removed by this technique.

\end{document}